\documentclass{amsart}
\usepackage{graphicx}
\usepackage{amsthm}
\usepackage{amsmath,latexsym}

\usepackage{amssymb}

\usepackage[latin1]{inputenc}
\usepackage{rotating}
\usepackage[a4paper, left=1.5cm, right=1.5cm, top=1.5cm, bottom=1.5cm]{geometry}
\usepackage{array}
\usepackage{multirow}

\usepackage[ruled,vlined]{algorithm2e}

\newtheorem{remark}{Remark}[section]
\newtheorem{corollary}{Corollary}[section]

\newtheorem{theo}{Theorem}[section]

\newcommand{\Z}{\mathbb{Z}}

\newcommand{\R}{\mathbb{R}}
\newcommand{\grob}{\mathcal{G}}

\newcommand{\dsum}{\displaystyle\sum}

\newcommand{\dprod}{\displaystyle\prod}

\begin{document}
\title{A Gröbner bases methodology for solving multiobjective polynomial integer programs}

\author{V. Blanco and J. Puerto}
\email{vblanco@us.es; puerto@us.es}
\address{Departamento de Estad\'istica e Investigaci\'on Operativa\\
Facultad de Matemáticas, Universidad de Sevilla, 41012 Sevilla, Spain.}
\thanks{This research was partially supported by
Spanish research grant numbers MTM2007- 67433-C02-01 and
P06-FQM-01366.}

\date{December 10, 2008}

\maketitle
\begin{abstract}
Multiobjective discrete programming is a well-known family of optimization
problems with a large spectrum of applications. The linear case has been tackled
by many authors during the last years. However, the polynomial case has not
been deeply studied due to its theoretical and computational difficulties. This
paper presents an algebraic approach for solving these problems. We propose a
methodology based on transforming the polynomial optimization problem in the
problem of solving one or more systems of polynomial equations and we use
certain Gröbner bases to solve these systems. Different transformations give
different methodologies that are analyzed and compared from a theoretical point
of view and by some computational experiments via the algorithms that they
induce.
\end{abstract}

\section{Introduction}
A multiobjective polynomial program consists of a finite set of
polynomial objective functions and a finite set of polynomial
constraints (in inequality or equation form), and solving that problem means obtaining the set of minimal elements in the feasible region defined by the constraints with respect to the partial order induced by the objective functions.

Polynomial programs have a wide spectrum of applications. Examples of them are capital budgeting~\cite{laughhunn70}, capacity
planning~\cite{shetty95}, optimization problems in graph
theory~\cite{beck00}, portfolio selection models with
discrete features~\cite{beasley03, jobst01} or chemical
engineering~\cite{ryoo-sahinidis95}, among many others. The reader is referred to \cite{li-sun06} for further applications.

Polynomial programming generalizes linear and quadratic programming and can serve as a tool to model engineering applications that
are expressed by polynomial equations. Even those problems with transcendental terms such as sin, log, and radicals can be reformulated by means of  Taylor series as a polynomial program. A survey of the publications on general nonlinear integer
programming can be found in \cite{cooper81}.

We study here multiobjective polynomial integer programs (MOPIP). Thus, we
assume that the feasible vectors have integer components and that there are more than one
objective function to be optimized. This change makes single-objective and multiobjective problems to be treated in a totally different manner, since the concept of solution is not the same.

In this paper, we introduce a new methodology for solving general
MOPIP based on the construction of reduced Gröbner bases of certain ideals
related to the problem and on solving triangular systems of
polynomial equations given by those bases. Gr\"obner bases were introduced  by Bruno Buchberger in 1965 in his PhD Thesis \cite{buchberger65}. He named it Gr\"obner basis paying tribute to his advisor Wolfgang
Gr\"obner. This theory emerged as a generalization, from the one
variable case to the multivariate polynomial case, of the
Euclidean algorithm, Gaussian elimination and the Sylvester
resultant. One of the outcomes of Gr\"obner
Bases Theory was its application to linear integer programming~\cite{conti-traverso91, hosten-sturmfels95, thomas97}. Later, Blanco and Puerto
 \cite{blanco-puerto07a} introduces a new notion of partial Gröbner basis for toric ideals in
 order to solve multiobjective linear integer programs. A different approach for solving linear integer programs was developed by Bertsimas et al.
\cite{bertsimas-perakis-tayur00} based on the application of Gröbner bases
for solving system of polynomial equations. This alternative use of
Gröbner bases is also used in the paper by Hägglof et al.
\cite{hagglof95} for solving continuous polynomial optimization
problems. Further details about Gröbner bases can be found in \cite{cox-little-oshea98,
cox-little-oshea92}.

We describe different approaches for solving MOPIP using Gröbner bases which are based on reducing the problem to several optimality conditions: the necessary Karush-Kuhn-Tucker, the Fritz-John and the multiobjective Fritz-John optimality conditions.

The paper is structured as follows. In the next section we give some preliminaries in multiobjective polynomial integer optimization. We present in Section \ref{sec:alg1} our first algorithm for solving MOPIP using only the triangularization property of lexicographic Gröbner bases. Section \ref{sec:alg2} is devoted to two different algorithms for solving MOPIP using a Chebyshev like scalarization and the Karush-Kuhn-Tucker or the Fritz-John optimality conditions. The last algorithm, based on the multiobjective Fritz-John optimality condition, is described in Section \ref{sec:alg4}. Finally, in Section \ref{sec:comp}, we compare the algorithms with the results of some computational
experiments and its analysis.

\section{The multiobjective integer polynomial problem}
\label{sec:mo}
The goal of this paper is to the solve multiobjective polynomial
integer programs (MOPIP):
\begin{equation}
\label{eq:mo-nonlinear}
\tag{$\mbox{\it MOPIP}_{\mathbf{f}, \mathbf{g}, \mathbf{h}}$}
\begin{array}{lrl}
\min & \quad (f_1(x), \ldots, f_k(x))& \\
s.t.& g_j(x) &\leq 0 \quad j =1, \ldots, m\\
 & h_r(x) &= 0 \quad r =1, \ldots, s\\
 & x &\in \Z_+^n
 \end{array}
 \end{equation}
with $f_1, \ldots, f_k, g_1, \ldots, g_m, h_1, \ldots, h_s$
polynomials in $\R[x_1, \ldots, x_n]$ and the constraints defining
a bounded feasible region. Therefore, from now on we deal with
$\mbox{\it MOPIP}_{\mathbf{f}, \mathbf{g}, \mathbf{h}}$ and we denote $\mathbf{f}=(f_1,
\ldots, f_k)$, $\mathbf{g}=(g_1, \ldots, g_m)$ and
$\mathbf{h}=(h_1, \ldots, h_r)$. If the problem had no equality (resp. inequality) constraints, we would denote it by $\mbox{\it MOPIP}_{\mathbf{f}, \mathbf{g}}$ (resp. $\mbox{\it MOPIP}_{\mathbf{f}, \mathbf{h}}$), avoiding the nonexistent term.

However, \ref{eq:mo-nonlinear} can be transformed to an
equivalent multiobjective polynomial binary problem. Since the feasible region $\{ x \in \R_+^n: g_j(x)\leq 0,\; h_r(x) =0,\,
j=1,\ldots, m, r=1, \ldots, s\}$ is assumed to be
bounded, it can be always embedded in a hypercube $\dprod_{i=1}^n
[0,u_i]^n$. Then, every component in $x$, $x_i$, has an
additional, but redundant, constraint $x_i \leq u_i$. We write
$x_i$ in binary form, introducing new binary variables $z_{ij}$ with values in $\{0,1\}$, $x_i = \sum_{j=0}^{\lfloor log u_i \rfloor} 2^j\,z_{ij}$, substituting every $x_i$ in \ref{eq:mo-nonlinear} we obtain an equivalent $0-1$ problem.

Then, from now on, without loss of generality, we restrict ourselves to multiobjective polynomial
binary programs ($MOPBP$) in the form:
\begin{equation}
\label{eq:mob-nonlinear}
\tag{$\mbox{\it MOPBP}_{\mathbf{f}, \mathbf{g}, \mathbf{h}}$}
\begin{array}{lrl}
\min & \quad (f_1(x), \ldots, f_k(x))& \\
s.t.& g_j(x) &\leq 0 \quad j =1, \ldots, m\\
  & h_r(x) &= 0 \quad r =1, \ldots, s\\
 & x &\in \{0,1\}^n
 \end{array}
 \end{equation}
 If the problem had no equality (resp. inequality) constraints, we would denote the problem by $\mbox{\it MOPBP}_{\mathbf{f}, \mathbf{g}}$ (resp. $\mbox{\it MOPBP}_{\mathbf{f}, \mathbf{h}}$), avoiding the nonexistent term.

The number of solutions of the above problem is finite, since the
decision space is finite. Thus, the number of feasible solutions
is, at most $|\{0,1\}^n| = 2^{n}$.

It is clear that $MOPBP_{\mathbf{f}, \mathbf{g}, \mathbf{h}}$ is
not a standard optimization problem since the objective function
is a $k$-coordinate vector, thus inducing a partial order among
its feasible solutions. Hence, solving the above problem requires
an alternative concept of solution, namely the set of nondominated
(or Pareto-optimal) points.

A feasible vector $\widehat{x} \in \R^n$ is said to be a
\textit{nondominated} (or Pareto optimal) solution of
$MOPIP_{{\mathbf{f}},{\mathbf{g}}}$ if there is no other feasible
vector $y$ such that
$$
f_j(y) \leq f_j(\widehat{x}) \qquad \forall j=1, \ldots, k
$$
with at least one strict inequality for some $j$. If $x$ is a
nondominated solution, the vector $\mathbf{f(x)} = (f_1(x),
\ldots, f_k(x))\in \R^k$ is called \textit{efficient}.

We say that a feasible solution, $y$, is dominated by a
feasible solution, $x$, if $f_i(x) \leq f_i(y)$ for all
$i=1,\ldots, k$ and $\mathbf{f}(x) \neq \mathbf{f}(y)$. We denote by $X_E$ the set of all nondominated
solutions for \ref{eq:mob-nonlinear} and by $Y_E$ the image under
the objective functions of $X_E$, that is, $Y_E = \{\mathbf{f}(x):
x \in X_E\}$. Note that $X_E$ is a subset of $\R^n$
(\emph{decision space}) and $Y_E$ is a subset of $\R^k$
(\emph{space of objectives}).

From the objective functions $\mathbf{f} = (f_1, \ldots, f_k)$, we
obtain a partial order on $\Z^n$ as follows:
$$
x \prec_\mathbf{f} y : \Longleftrightarrow \mathbf{f}(x)
\lvertneqq \mathbf{f}(y) \mbox{ or } x=y.
$$
Note that since $\mathbf{f}: \R^n \rightarrow \R^k$, the above
relation is not complete. Hence, there may exist incomparable
vectors.

In the following sections we describe some algorithms for solving MOPIP using tools from algebraic geometry. In particular, in each of these methods, we transform our problem in a certain system of polynomial equations, and we use Gröbner bases to solve it.

\section{Obtaining nondominated solutions solving systems of polynomial equations}
\label{sec:alg1}
In this section we present the first
approach for solving multiobjective polynomial integer programs
using Gröbner bases. For this method, we transform the program in a system of polynomial equations that encodes the set of feasible solutions and its objective values. Solving that system in the objective values, and then, selecting the minimal ones in the partial componentwise order, allows us to obtain the associate feasible vectors, thus, the nondominated solutions.

Through this section we solve $\mbox{\it MOPBP}_{\mathbf{f}, \mathbf{h}}$. Without loss of generality, we reduce the general problem to the problem without inequality constraints since we can transform inequality constraints to equality constrains as follows:
\begin{equation}
\label{eq:slacks}
g(x) \leq 0 \Longleftrightarrow g(x)+ z^2 =0, z\in \R.
\end{equation}
where the quadratic term, $z^2$, assures the nonnegativity of the slack
variable and then, less than or equal to type inequality. Initially, we suppose that all the variables are binary. In Remark \ref{remark1} we describe how to modify the algorithm to incorporate the above slack variables.

This approach consists of transforming $\mbox{\it MOPBP}_{\mathbf{f}, \mathbf{h}}$ to an equivalent problem such that the
objective functions are part of the constraints. For this
transformation, we add $k$ new variables, $y_1, \ldots, y_k$ to the
problem, encoding the objective values for all feasible solutions. The modified problem is:
\begin{equation}
\label{eq:nonlinear-alg1_2}
\begin{array}{lrl}
\min & \quad (y_1, \ldots, y_k)& \\
s.t. & h_r(x) &= 0 \quad r =1, \ldots, s\\
 & y_j - f_j(x) & = 0 \quad j =1, \ldots, k\\
 & x_i(x_i -1) & = 0 \quad i =1, \ldots, n\\
 & y \in\R^k& x \in\R^n\\
 \end{array}
 \end{equation}
where integrality constraints are codified as quadratic constraints so, $\mbox{\it MOPBP}_{\mathbf{f}, \mathbf{h}}$ is a polynomial continuous problem.

 The algorithm consists of, first, obtaining the set of feasible
 solutions of Problem \eqref{eq:nonlinear-alg1_2} in the $y$
 variables; then, selecting from that set those solutions that are
 minimal with respect to the componentwise order, obtaining the
 set of efficient solutions of $\mbox{\it MOPBP}_{\mathbf{f}, \mathbf{h}}$. The feasible solutions in the
 $x$-variables associated to those efficient solutions correspond
 with the nondominated solutions of $\mbox{\it MOPBP}_{\mathbf{f}, \mathbf{h}}$.

 Then, first, we concentrate in describing a procedure for solving the
 system of polynomial equations that encodes the feasible region
 of Problem \eqref{eq:nonlinear-alg1_2}, i.e. the solutions of
\begin{equation}
\label{eq:system}
 \begin{array}{lll}
h_r(x) &= 0 & \mbox{ for all } r=1,\ldots, s\\
y_j - f_j(x) &= 0 & \mbox{ for all }  j=1,\ldots, k\\
x_i(x_i-1) &=0 & \mbox{ for all }  i=1,\ldots, n.
 \end{array}
\end{equation}
For analyzing the system \eqref{eq:system} we use Gröbner bases as a tool for solving systems of polynomial equations. Further details can be found in the book by Sturmfels~\cite{sturmfels02}.

The set of solutions of \eqref{eq:system} coincides with the affine variety of the following polynomial ideal in
$\R[y_1, \ldots, y_k, x_1, \ldots, x_n]$:
$$
I = \langle h_1(x), \ldots, h_m(x), y_1 - f_1(x), \ldots, y_k -
f_k(x), x_1(x_1-1), \ldots, x_n(x_n-1) \rangle.
$$
Note that $I$ is a zero dimensional ideal since the number
of solutions of the equations that define $I$ is finite. Let $V(I)$ denote the affine variety of $I$. If we restrict $I$ to the family of variables $x$ (resp. $y$) the variety $V(I\cap \R[x_1, \ldots, x_n])$ (resp. $V(I\cap
\R[y_1, \ldots, y_k])$) encodes the set of feasible solutions (resp. the set of possible objective values) for that problem.

Applying the elimination property, the
reduced Gröbner basis for $I$, $\grob$, with respect to the
lexicographical ordering with $y_k \prec \cdots \prec y_1 \prec
x_n \prec \cdots \prec x_1$ gives us a method for solving system \eqref{eq:system} sequantially, i.e., solving in one indeterminate at a time. Explicitly, the shape of
$\grob$ is:
\begin{enumerate}
\item[$1)$] $\grob$ contains one polynomial in $\R[y_k]$: $p_k(y_k)$
\item[$2)$] $\grob$ contains one or several polynomials in $\R[y_{k-1},
y_k]$ : $p_{k-1}^1(y_{k-1},y_k), \ldots, p_{k-1}^{m_{k-1}}(y_{k-1},y_k)$.
\item[$\vdots$]
\item[$k+1)$] $\grob$ contains one or several polynomials in $\R[x_n, y_{1},
\ldots, y_k]$ : $q_{n}^1(x_n,\mathbf{y}), \ldots, q_{n}^{s_n}(x_n,\mathbf{y})$.
\item[$\vdots$]
\item[$k+n)$] $\grob$ contains one or several polynomials in $\R[x_n, y_{1},
\ldots, y_k]$ : $q_{1}^1(x_1, \ldots, x_n,\mathbf{y}), \ldots, q_{n}^{s_1}(x_1, \ldots, x_n,\mathbf{y})$.
\end{enumerate}
Then, with this structure of $\grob$, we can solve, in a first step, the
system in the $y$ variables using those polynomial in $\grob$ that
only involve this family of variables as follows: we first solve
for $y_k$ in $p_k(y_k)=0$, obtaining the solutions: $y_k^1, y_k^2,
\ldots$. Then, for fixed $y_k^r$, we find the common roots of
$p_{k-1}^1, p_{k-1}^2, \ldots$ getting solutions $y_{k-1, r}^1,
y_{k-1, r}^2, \ldots$ and so on, until we have obtained the roots for
$p_1(y_1, \ldots, y_k)$. Note that at each step we only solve
one-variable polynomial equations.

We denote by $\Omega$ the above set of solutions in vector form
\begin{eqnarray*}
\begin{split}
\Omega = \{ (\hat{y}_1, \ldots, \hat{y}_k) : &p_k(\hat{y}_k)=0, p_{k-1}^1(\hat{y}_{k-1}, \hat{y}_k)=0, \ldots,
p_{k-1}^{m_{k-1}}(\hat{y}_{k-1}, \hat{y}_k)=0, \ldots\\
& p^1_1(\hat{y}_1, \hat{y}_2, \ldots, \hat{y}_k)=0, \ldots,
p^{m_1}_1(\hat{y}_1, \hat{y}_2, \ldots, \hat{y}_k)=0 \}.
\end{split}
\end{eqnarray*}
As we stated above, $\Omega$ is the set of all possible values of the objective functions at the feasible solutions of $\mbox{\it MOPBP}_{\mathbf{f}, \mathbf{h}}$. We are looking for the nondominated
solutions that are associated with the efficient solutions. From
$\Omega$, we can select the efficient solutions as those that are
minimal with respect to the componentwise order in $\R^k$. So, we
can extract from $\Omega$ the set of efficient solutions, $Y_E$:
$$
Y_E = \{(y_1^*, \ldots, y_k^*)  \in \Omega: \not\exists (y_1^\prime, \ldots, y_k^\prime)  \in \Omega \text{ with } y_j^\prime\leq y_j^* \text{ for } j=1, \ldots, k \text{ and } (y_1^\prime, \ldots, y_k^\prime) \neq (y_1^*, \ldots, y_k^*)\}
$$
Once we have obtained the solutions in the $y$ variables that are
efficient solutions for $\mbox{\it MOPBP}_{\mathbf{f}, \mathbf{h}}$, we compute with
an analogous procedure the nondominated solutions associated to
the $y$-values in $Y_E$. It consists of solving the triangular system given by
$\grob$ for the polynomial where the $x$-variables appear once the values for the $y$-variables are fixed to be each of the
vectors in $Y_E$.

A pseudocode for this procedure is described in Algorithm
\ref{alg1}.

\begin{algorithm}[h]
\label{alg1} \SetLine \SetKwInOut{Input}{Input}
\SetKwInOut{Output}{Output}

\Input{$f_1, \ldots, f_k$, $h_1, \ldots h_s \in \R[x_1, \ldots,
x_n]$}

\textbf{Initialization: } $I= \langle f_1-y_1, \ldots, f_k-y_k,
h_1, \ldots, h_s, x_1(x_1-1), \ldots, x_n(x_n-1)\rangle$.\\

\textbf{Algorithm: }

\begin{description}
\item[Step 1.] Compute a Gröbner basis, $G$, for $I$ with respect to a lexicographic order with $y_k \prec \cdots \prec y_1 \prec x_n
\prec \cdots \prec x_1$.
\item[Step 2.] Let $G^y_l = G \cap \R[y_{l+1}, \ldots, y_k]$ be a
Gröbner basis for $I^y_l=I \cap \R[y_{l+1}, \ldots, y_k]$, for
$l=0, \ldots, k-1$. (By the Elimination Property).
\begin{enumerate}
\item Find all $\hat{y}_k \in V(G_{k-1}^y)$.
\item Extend every $\hat{y}_k$ to $(\hat{y}_{k-1}, \hat{y}_k) \in
V(G_{k-2}^y)$.
\item[$\vdots$]
\item[($k-1$)] Extend every $(\hat{y}_3, \ldots, \hat{y}_k)$ to $(\hat{y}_2, \hat{y}_3, \ldots, \hat{y}_k) \in
V(G_{1}^y)$.
\item[($k$)] Find all $\hat{y}_1$ such that $(\hat{y}_1, \ldots, \hat{y}_k) \in
V(G_{0}^y)$.
\end{enumerate}
\item[Step 3.] Select from $V(G_{0}^y)$ the minimal vectors with
respect to the usual componentwise order in $\R^k$. Set $Y_E$ this subset.
\item[Step 4.] Let $G_l = G \cap \R[y_{1}, \ldots, y_k, x_{l+1}, \ldots, x_n]$ a
Gröbner basis for $I_l \cap  \R[y_{1}, \ldots, y_k, x_{l+1},
\ldots, x_n]$, for $l=0, \ldots, n-1$. (By the Elimination
Property). Denote by $S_l = \{(\hat{y}_1, \ldots, \hat{y}_k, \hat{x}_{l+1},\ldots, \hat{x}_n): (\hat{y}_1, \ldots, \hat{y}_k) \in
Y_E, \text{ and } \exists (x_1, \ldots, x_l) \text{ such that } (x_1, \ldots, x_n) \text{ is feasible }\}$ for $l=0, \ldots, n-1$.
\begin{enumerate}
\item Find all $\hat{x}_n$ such that $(\hat{y}_1, \ldots, \hat{y}_k, \hat{x}_n) \in V(G_{n-1}) \cap S_{n-1}$.
\item Extend every $\hat{x}_n$ to $((\hat{y}_1, \ldots, \hat{y}_k, \hat{x}_{n-1}, \hat{x}_n) \in
V(G_{n-2})\cap S_{n-2}$.
\item[$\vdots$]
\item[($n-1$)] Extend every $(\hat{y}_1, \ldots, \hat{y}_k, \hat{x}_3, \ldots, \hat{x}_n)$ to $(\hat{y}_1, \ldots, \hat{y}_k, \hat{x}_2, \hat{x}_3,\ldots, x_n) \in
V(G_{1})\cap S_{1}$.
\item[($n$)] Find all $\hat{x}_1$ such that $(\hat{y}_1, \ldots, \hat{y}_k, \hat{x}_1, \ldots, \hat{x}_n) \in
V(G_{0})\cap S_{0}$.
\end{enumerate}
Set $X_E = \pi_x(V(G_{0})\cap S_{0})$, where $\pi_x$ denotes the
projection over the $x$-variables.
\end{description}
\Output{$Y_E$ the set of efficient solutions and $X_E$ the set of
nondominated solutions for $MOPBP_{\mathbf{f}, \mathbf{h}}$.}
 \caption{Solving MOPIP by solving systems of polynomial equations}
\end{algorithm}
\begin{theo}
Algorithm \ref{alg1} either provides all nondominated and
efficient solutions or provides a certificate of infeasibility
whenever $G=\{1\}$.
\end{theo}

\begin{proof}
Suppose that $G\neq \{1\}$. Then, $G_{k-1}^y$ has exactly one element, namely $p(y_k)$. This follows from the observation that $I \cap
\R[y_k]$ is a polynomial ideal in one variable, and therefore,
needs only one generator.

Solving $p(y_k)=0$ we obtain every $\hat{y}_k \in V(G_{k-1}^y)$.
Sequentially we obtain $\hat{y}_{k-1}$ extending $\hat{y}_{k}$ to the partial
solutions $(\hat{y}_{k-1}, \hat{y}_k)$ in $V(G_{k-1}^y)$ and so
on.

By the Extension Theorem, this is always possible in our case.

Continuing in this way and applying the Extension Theorem, we can
obtain all solutions $(\hat{y}_{1}, \ldots, \hat{y}_k)$ in $V(G
\cap \R[y_1, \ldots, y_k]$. These vectors are all the possible
objective values for all feasible solutions of the problem.
Selecting from $V(G \cap \R[y_1, \ldots, y_k])$ those solutions
that are not dominated in the componentwise order in $\R^k$, we
obtain $Y_E$.

Following a similar scheme in the $x$- variables, we have
the set $V(G_{0})\cap S_{0}^*$ encoding all efficient (in the
first $k$ coordinates) and nondominated (in the last $n$
coordinates) solutions.

Finally, if $G=\{1\}$, then, the ideal $I$ coincides with
$\R[y_1, \ldots, y_k, x_1, \ldots, x_n]$, indicating that $V(I)$
is empty (it is the set of the common roots of all polynomials in
$\R[y_1, \ldots, y_k, x_1, \ldots, x_n]$). Then, we have an
infeasible integer problem.
\end{proof}

\begin{remark}
\label{remark1}
In the case when we have added slack variables, as explained in \eqref{eq:slacks}, we slightly modify the above algorithm solving first in the slack variables and selecting those solutions that are real numbers. Then continue with the procedure described in Algorithm \ref{alg1}.
\end{remark}

\begin{remark}
\label{remark2}The Gröbner basis, $\grob$, computed for solving the system of polynomial equations can be computed with respect to any other elimination ordering. The only conditions that are required for that ordering is that it allows to separate the family of $x-$variables from the family of $y$-variables and such that the system of polynomials given by that basis allows solving first for the $y$-variables and then for the $x$-variables sequentially.
\end{remark}

\section{Obtaining nondominated solutions by the Chebyshev norm
approach} \label{sec:alg2}

In this section we describe two more methods for solving MOPIP based on a different rationale, namely scalarizing the multiobjective problem and solving it as a parametric single-objective problem. We propose a methodology based on the application of optimality conditions to a family of
single-objective problems related to our original multiobjective
problem. The methods consist of two main steps: a first step where the multiobjective problem is
scalarized to a family of single-objective problems such that each
nondominated solution is an optimal solution for at least one of
the single-objective problems in that family; and a second step that
consists of applying necessary optimality conditions to each one of the
problems in the family, to obtain their optimal solutions. Those solutions are only candidates to be nondominated
solutions of the multiobjective problem since we just use necessary
conditions.

For the first step, the scalarization, we use a weighted Chebyshev
norm approach. Other weighted sum approaches could be used to
transform the multiobjective problem in a family of single-objective
problems whose set of solutions contains the set of nondominated solutions of our problem. However, the
Chebyshev approach seems to be rather adequate since it does not require to impose extra hypothesis to the problem. This
approach can be improved for problems satisfying convexity
conditions, where alternative well-known results can be applied (see
\cite{jahn04} for further details).

For the second step, we use the Fritz-John and Karush-Kuhn-Tucker necessary optimality
conditions, giving us two different approaches. In this section we describe both
methodologies since each of them has its own advantages over the
other.

For applying the Chebyshev norm scalarization, we use the
following result that states how to transform our problem to a
family of single-objective problems, and how to obtain
nondominated solutions from the optimal solution of those
single-objective problems. Further details and proofs of this
result can be found in \cite{jahn04}.

\begin{theo}[Corollary 11.21 in \cite{jahn04}] Let \ref{eq:mob-nonlinear} be feasible. $x^*$ is a nondominated
solution of \ref{eq:mob-nonlinear} if and only if there
are positive real numbers $\omega_1, \ldots, \omega_k
>0$ such that $x^*$ is an image unique solution of the following weighted
Chebyshev approximation problem:
\begin{equation}
\label{eq:wmo-nonlinear} \tag{$P_\omega$}
\begin{array}{lrl}
\min & \gamma & \\
s.t.& \omega_i\,(f_i(x)-\hat{y}_i) - \gamma &\leq 0 \quad i =1, \ldots, k\\
 & g_j(x) &\leq 0 \quad j =1, \ldots, m\\
 & h_r(x) &= 0 \quad r =1, \ldots, s\\
 & x_i(x_i-1) &=0  \quad i =1, \ldots, n\\
 & \gamma \in \R & x \in \R^n
 \end{array}
 \end{equation}
 where $\hat{y}=(\hat{y}_1, \ldots, \hat{y}_k)\in \R^k$ is a lower bound of $f=(f_1, \ldots,
 f_k)$, i.e., $\hat{y}_i \leq f_i(x)$  for all feasible solution $x$ and $i=1, \ldots,
 k$.
\end{theo}

According to the above result, every nondominated solution of
\ref{eq:mob-nonlinear} is the unique solution of
\eqref{eq:wmo-nonlinear} for some $\mathbf{\omega} >0$. We apply,
in the second step, necessary optimality conditions for obtaining
the optimal solutions for those problems (taking $\omega$ as
parameters). These solutions are candidates to be nondominated
solutions of our original problem. Actually, every nondominated
solution is among those candidates.

In the following subsections we describe the above-mentioned two
methodologies for obtaining the optimal solutions for the
scalarized problems \eqref{eq:wmo-nonlinear} for each $\omega$.

\subsection{The Chebyshev-Karush-Kuhn-Tucker approach:}

The first optimality conditions that we apply are the
Karush-Kuhn-Tucker (KKT) necessary optimality conditions, that were
stated, for the general case, as follows (see e.g. \cite{bazaraa93} for further details):

\begin{theo}[KKT necessary conditions]
\label{kkt_cond} Consider the problem:
\begin{equation}
\label{eq:nl-kkt}
\begin{array}{lrl}
\min & f(x) & \\
s.t.& g_j(x) &\leq 0 \quad  =1, \ldots, m\\
 & h_r(x) &= 0 \quad r=1, \ldots, s\\
 & x &\in \R^n
 \end{array}
 \end{equation}
Let $x^*$ be a feasible solution, and let $J=\{j: g_j(x^*)=0\}$.
Suppose that $f$ and $g_j$, for $j=1, \ldots, m$, are
differentiable at $x^*$,  that $g_j$, for $j\not\in J$, is
continuous at $x^*$,and that $h_r$, for $r=1, \ldots, s$, is
continuously differentiable at $x^*$. Further suppose that $\nabla
g_j$, for $j\in I$, and $\nabla h_r$, for $r=1, \ldots, s$, are
linearly independent (regularity conditions). If $x^*$ solves
Problem \ref{eq:nl-kkt} locally, then there exist scalars
$\lambda_j$, for $j=1, \ldots, m$, and $\mu_r$, for $r=1, \ldots,
s$, such that
\begin{equation}
\label{eq:kkt} \tag{KKT}
\begin{array}{rll}
 \nabla f(x^*) + \dsum_{j=1}^m
\lambda_j\,\nabla g_j(x^*)\; +& \dsum_{r=1}^s \mu_r\,\nabla
h_r(x^*)
&=0\\
\lambda_j\,g_j(x^*)\;=&0 & \mbox{for $j=1, \ldots, m$}\\
\lambda_j \;\geq& 0 & \mbox{for $j=1, \ldots, m$}
\end{array}
\end{equation}
\end{theo}

From the above theorem the
candidates to be optimal solutions for Problem \eqref{eq:nl-kkt}
are those that either satisfy the KKT conditions (in the case where all the functions
involved in Problem \eqref{eq:nl-kkt} are polynomials, this is a system of polynomial equations) or do not satisfy
the regularity conditions. Note that these two sets are, in general, not disjoint.

Regularity conditions can also be formulated as a system of polynomial equations when the involved functions are all polynomials. Let
$x^*$ be a feasible solution for Problem \eqref{eq:nl-kkt}, $x^*$
does not verify the regularity conditions if there exist scalars
$\lambda_j$, for $j\in J$, and $\mu_r$, for $r=1, \ldots, s$, not
all equal to zero, such that:
\begin{equation}
\label{regularity}\tag{Non-Regularity}
 \dsum_{j\in I} \lambda_j\nabla g_j + \dsum_{r=1}^s \mu_r\nabla
h_r = 0
\end{equation}

The above discussion justifies the following result.
\begin{corollary}
Let $x^*$ be a nondominated solution for \ref{eq:mob-nonlinear}. Then, $x^*$ is a solution of the systems of polynomial equations \eqref{kkt:pw} or \eqref{nr:pw}, for some $\omega>0$.

\begin{eqnarray}
\label{kkt:pw}
\left.
\begin{split}
1 - \dsum_{i=1}^k \nu_i = 0&\\
\dsum_{i=1}^k \nu_i \omega_i \,\nabla f_i(x) + \dsum_{j=1}^m \lambda_j \nabla g_j(x) + \dsum_{r=1}^s \mu_j \nabla h_r(x)+ \dsum_{l=1}^n \beta_l\,\delta_{il}(2x_i-1) = 0&\\
\nu_i\,\omega_i\,(f_i(x)-\hat{y}_i) - \gamma =0, & \; i=1, \ldots, k\\
\end{split}\right\}
\end{eqnarray}
 for some $\lambda_j\,g_j(x) =0$, for $j=1, \ldots, m$, $\lambda_j \geq 0$, for $j=1, \ldots, m$, and $\nu_i \geq 0$, for $i=1, \ldots, k$.

\begin{eqnarray}
\label{nr:pw}
\left.
\begin{split}
\dsum_{i=1}^k \nu_i = 0 &\\
 \dsum_{i=1}^k \nu_i \omega_i \,\nabla f_i(x) + \dsum_{j=1}^m
\lambda_j \nabla g_j(x) + \dsum_{r=1}^s
\mu_j \nabla h_r(x)+ \dsum_{l=1}^n \beta_l\,\delta_{il}(2x_i-1) = 0&\\
\omega_i\,(f_i(x)-\hat{y}_i) - \gamma \leq 0, & \; i=1, \ldots,
k\\
\end{split}\right\}
\end{eqnarray}
with $x\in \R^n$ such that $g_j(x) \leq 0$, for $j=1, \ldots, m$, $h_r(x) = 0$, for $r=1, \ldots, s$\\
$\lambda_j \geq 0$, for $j=1, \ldots, m$, and $\nu_i \geq 0$, for $i=1, \ldots, k$.
\end{corollary}

Let $X_E^{KKT}$ denote  the set of solutions, in the
$x$-variables, of system \eqref{kkt:pw} and let $X_E^{NR}$ denote the set of solutions, in the
$x$-variables, of system \eqref{nr:pw} (the problem is solved
avoiding inequality constraints, then every solution is evaluated to
check if it satisfies the inequality constraints).

For solving these systems (Chebyshev-KKT and Non-Regulariry), we use a
Gröbner bases approach. Let $I$ be the ideal generated by the
involved equations.

Let us consider a lexicographical order over the monomials in $\R[\mathbf{x}, \gamma, \mathbf{\lambda}, \mathbf{\nu}, \mathbf{\mu}, \mathbf{\beta}]$
such that $\mathbf{x} \prec \gamma \prec \mathbf{\lambda} \prec
\mathbf{\nu} \prec \mathbf{\mu} \prec \mathbf{\beta}$. Then, the
Gröbner basis, $\grob$, for $I$ with this order has the following
triangular shape:
\begin{itemize}
\item $\grob$ contains one polynomial in $\R[x_n]$: $p_n(x_n)$
\item $\grob$ contains one or several polynomials in $\R[x_{n-1},
x_n]$ : $p_{n-1}^1(x_{n-1}, x_n), \ldots, p_{n-1}^{m_1}(x_{n-1}, x_n)$.
\item[$\vdots$]
\item $\grob$ contains one or several polynomials in $\R[\mathbf{x}]$ : $p_{1}^1(x_1, \ldots, x_n), \ldots, p_{1}^{m_n}(x_{1}, \ldots, x_n)$.
\item The remaining polynomials involve variables $\mathbf{x}$ and at least one $\gamma$,
$\lambda$, $\mu$, $\nu$ or $\beta$.
\end{itemize}

We are interested in finding only the values for the
$x$-variables, so, we avoid the polynomials in $\grob$ that
involve any of the other auxiliary variables. In general, we are not able to discuss about the values of the parameters $\gamma$,
$\lambda$, $\mu$, $\nu$ and $\beta$. Needless to say that in those cases when we can do it, some values of $\mathbf{x}$ may be discarded simplifying the process. We denote by
$\grob^x$ the subset of $\grob$ that contains only polynomials in
the $x$-variables. By the Extension Theorem, $\grob^x$ is a
Gröbner basis for $I \cap \R[x_1, \ldots , x_n]$.

Solving the system given by $\grob^x$ and checking feasibility of those solutions, we obtain as solutions those of our KKT or Non-Regularity original systems.

It is clear that the set of nondominated solutions of our
problem is a subset of $X_E^{KKT} \cup X_E^{NR}$, since either a
solution is regular, and then, KKT conditions are applicable or it
satisfies the non regularity conditions. However, the set
$X_E^{KKT} \cup X_E^{NR}$ may contain dominated solutions, so, at the end we must remove the dominated ones to get only $X_E$.

The steps to solve Problem \ref{eq:mob-nonlinear} using the Chebyshev-KKT
approach are summarized in Algorithm \ref{alg:kkt}.
\begin{algorithm}[h]
\label{alg:kkt} \SetLine \SetKwInOut{Input}{Input}
\SetKwInOut{Output}{Output}

\Input{$f_1, \ldots, f_k$, $g_1, \ldots g_m, h_1, \ldots, h_r \in \R[x_1, \ldots,
x_n]$}

\textbf{Algorithm: }
\begin{description}
\item[Step 1] Formulate the Chebyshev scalarization of \ref{eq:mob-nonlinear}. (Problem \eqref{eq:wmo-nonlinear})
\item[Step 2] Solve System \eqref{kkt:pw} in the $x$-variables: $X_E^{KKT}$.
\item[Step 3] Solve System \eqref{nr:pw} in the $x$-variables: $X_E^{NR}$.
\item[Step 4] Remove from $X_E^{KKT} \cup X_E^{NR}$ the subset of dominated solutions: $X_E$.
\end{description}
\Output{$X_E$ the set of nondominated solutions for \ref{eq:mob-nonlinear}}
 \caption{Summary of the procedure for solving MOPBP using Chebyshev-KKT approach.}
\end{algorithm}

\begin{theo}
Algorithm \ref{alg:kkt} solves Problem \ref{eq:mob-nonlinear} in a finite number of steps.
\end{theo}

\subsection{The Chebyshev-Fritz-John approach}
Analogously to the previous approach, once we have
scalarized the original multiobjective problem to a family of
single-objective problems, in this section we apply the Fritz-John (FJ) conditions to all the
problems in this family. The following well-known result justifies the use of FJ conditions to obtain candidates to optimal
solutions for single-objective problems. Proofs and further
details can be found in \cite{bazaraa93}.

\begin{theo}[FJ necessary conditions]
\label{fj_cond} Consider the problem:
\begin{equation}
\label{eq:nl-fj}
\begin{array}{lrl}
\min & f(x) & \\
s.t.& g_j(x) &\leq 0 \quad  =1, \ldots, m\\
 & h_r(x) &= 0 \quad r=1, \ldots, s\\
 & x &\in \R^n
 \end{array}
 \end{equation}
Let $x^*$ be a feasible solution, and let $J=\{j: g_j(x^*)=0\}$.
Suppose that $f$ and $g_j$, for $j=1, \ldots, m$, are
differentiable at $x^*$, and that $h_r$, for $r=1, \ldots, s$, is
continuously differentiable at $x^*$. If $x^*$ locally solves
Problem \eqref{eq:nl-fj}, then there exist scalars
$\lambda_j$, for $j=1, \ldots, m$, and $\mu_r$, for $r=1, \ldots,
s$, such that
\begin{equation}
\label{eq:fj} \tag{FJ}
\begin{array}{rll}
\lambda_0 \nabla f(x^*) + \dsum_{j=1}^m \lambda_j\,\nabla
g_j(x^*)\; +& \dsum_{r=1}^s \mu_r\,\nabla h_r(x^*)
&=0\\
\lambda_j\,g_j(x^*)\;=&0 & \mbox{for $j=1, \ldots, m$}\\
\lambda_j \;\geq& 0 & \mbox{for $j=1, \ldots, m$}\\
(\lambda_0, \mathbf{\lambda}, \mathbf{\mu}) &\neq (0, \mathbf{0},
\mathbf{0})
\end{array}
\end{equation}
\end{theo}

Note that, in the FJ conditions, regularity conditions are not
required to set the result.

\begin{corollary}
Let $x^*$ be a nondominated solution for \ref{eq:mob-nonlinear}. Then, $x^*$ is a solution of the system of polynomial equations \eqref{fj:pw} for some $\nu_i, \lambda_j, \mu_r, \beta_l$, for $i=1, \ldots, k, l=1, \ldots,
n$, $j=1, \ldots, m$, $r=1, \ldots, s$ and $\omega>0$.
\begin{eqnarray}
\label{fj:pw}
\left.
\begin{split}
\lambda_0 - \dsum_{i=1}^k \nu_i = 0\\
 \dsum_{i=1}^k \nu_i \omega_i \,\nabla f_i(x) + \dsum_{j=1}^m
\lambda_j \nabla g_j(x) + \dsum_{r=1}^s
\mu_j \nabla h_r(x)+ \dsum_{l=1}^n \beta_l\,\delta_{il}(2x_i-1) = 0\\
\nu_i\,(\omega_i\,(f_i(x)-\hat{y}_i) - \gamma) =0, & i=1, \ldots,
k\\
\lambda_j\,g_j(x) =0, & \; j=0, \ldots, m
\end{split}\right\}
\end{eqnarray}
where $\lambda_j \geq 0$, for $j=1, \ldots, m$, $\nu_i \geq 0$, for $i=1, \ldots, k$ and $\delta_{li}$,
for $l,i = 1, \ldots, n$, denotes the Kronecker delta function, $\delta_{li}=\left\{\begin{array}{ll}
1 & \mbox{if $l=i$}\\
0 & \mbox{otherwise}
\end{array}\right.$.
\end{corollary}

Let $X_E^{FJ}$ denote the set of solutions, in the
$x$-variables, that are feasible solutions of \ref{eq:mob-nonlinear} and solutions of system \eqref{fj:pw}.

The set of nondominated solutions of our
problem is a subset of $X_E^{FJ}$, since every nondominated
solution is an optimal solution for some problem in the form of \eqref{fj:pw}, and
every solution of this single-objective problem is a solution of
the FJ system.

However, dominated solutions may appear in the set of solutions of
\eqref{fj:pw}, so, a final elimination process is to be performed to select only the
nondominated solutions.

The steps to solve \ref{eq:mob-nonlinear} using the Chebyshev-FJ
approach are summarized in Algorithm \ref{alg:fj}.

\begin{algorithm}[h]
\label{alg:fj} \SetLine \SetKwInOut{Input}{Input}
\SetKwInOut{Output}{Output}

\Input{$f_1, \ldots, f_k$, $g_1, \ldots g_m, h_1, \ldots, h_r \in \R[x_1, \ldots,
x_n]$}

\textbf{Algorithm: }
\begin{description}
\item[Step 1] Formulate the Chebyshev scalarization of \ref{eq:mob-nonlinear}. (Problem \eqref{eq:wmo-nonlinear})
\item[Step 2] Solve system \eqref{fj:pw} in the $x$-variables for any value of $\omega > 0$: $X_E^{FJ}$.
\item[Step 3] Remove from $X_E^{FJ}$ the set of dominated solutions: $X_E$.
\end{description}

\Output{$X_E$ the set of nondominated solutions for Problem
\eqref{eq:mob-nonlinear}}

 \caption{Summary of the procedure for solving MOPBP using the Chebyshev-FJ approach.}
\end{algorithm}

\begin{theo}
Algorithm \ref{alg:fj} solves \ref{eq:mob-nonlinear} in a finite number of steps.
\end{theo}

The last part of the section is devoted to show how to
solve the Chebyshev-FJ system using Gröbner bases.

Consider the following polynomial ideal

{\small
$I = \langle \lambda_0 -
\dsum_{i=1}^k \nu_i, \dsum_{i=1}^k \nu_i \omega_i \,\nabla f_i(x)
+ \dsum_{j=1}^m \lambda_j \nabla g_j(x) + \dsum_{r=1}^s \mu_j
\nabla h_r(x)+ \dsum_{l=1}^n \beta_l\,\delta_{il}(2x_i-1),
\nu_1\,(\omega_1\,(f_1(x)-\hat{y}_1) - \gamma),\ldots,
\nu_k\,(\omega_k\,(f_k(x)-\hat{y}_k) - \gamma), \lambda_1\,g_1(x),
\ldots, \lambda_m\,g_m(x)\rangle $}

in the polynomial ring $\R[\mathbf{x}, \gamma, \mathbf{\lambda}, \mathbf{\nu}, \mathbf{\mu}, \mathbf{\beta}]$.

Let us consider a lexicographical order over the monomials in $\R[\mathbf{x}, \gamma, \mathbf{\lambda}, \mathbf{\nu}, \mathbf{\mu}, \mathbf{\beta}]$ such that $\mathbf{x} \prec \gamma \prec \mathbf{\lambda} \prec
\mathbf{\nu} \prec \mathbf{\mu} \prec \mathbf{\beta}$. Then, the
Gröbner basis, $\grob$, for $I$ with this order has the following
triangular shape:
\begin{itemize}
\item $\grob$ contains one polynomial in $\R[x_n]$: $p_n(x_n)$
\item $\grob$ contains one or several polynomials in $\R[x_{n-1},
x_n]$ : $p_{n-1}^1(x_{n-1}, x_n), \ldots, p_{n-1}^{m_1}(x_{n-1}, x_n)$
\item[$\cdots$]
\item $\grob$ contains one or several polynomials in $\R[\mathbf{x}]$ : $p_{1}^1(x_1, \ldots, x_n), \ldots, p_{1}^{m_n}(x_{1}, \ldots, x_n)$
\item The remainder polynomials involve variables $\mathbf{x}$ and at least one of $\gamma$,
$\lambda$, $\mu$, $\nu$ or $\beta$.
\end{itemize}

We are interested in finding only the values for the
$x$-variables, so, we avoid the polynomials in $\grob$ that
involve any of the other auxiliary variables. We denote by
$\grob^x$ the subset of $\grob$ that contains only all the polynomials in
the $x$-variables. By the Extension Theorem, $\grob^x$ is a
Gröbner basis for $I \cap \R[x_1, \ldots , x_n]$.

Solving the system given by $\grob^x$, we obtain as solutions,
those of our FJ original system.

\begin{remark}[Convex Case]
In the special case where both objective functions and constraints are convex, sufficient KKT conditions can be
applied. If the feasible solution $x^*$ satisfies KKT
conditions, and all objective and constraints functions are
convex, then $x^*$ is a nondominated solution. As a particular
case, this situation is applicable to linear problems.

In this case, we may choose a linear scalarization instead of the Chebyshev scalarization. With this alternative approach, the scalarized problem is
\begin{equation*}
\label{eq:linmo-nonlinear}
\begin{array}{lrl}
\min & \dsum_{s=1}^k t_s\,f_s(x) & \\
s.t.& g_j(x) &\leq 0 \quad j =1, \ldots, m\\
 & h_r(x) &= 0 \quad r =1, \ldots, s\\
 & x_i(x_i-1) &=0  \quad i =1, \ldots, n\\
 \end{array}
\end{equation*}
for $t_1, \ldots, t_k >0$.

Then, by Corollary 11.19 in \cite{jahn04}, and denoting by $S$ the feasible region, if $f(S)+\R^k_+$ is convex, then each $x^*$ is a nondominated solution if and only if $x^*$ is a solution of Problem \ref{eq:linmo-nonlinear} for some $t_1, \ldots, t_k>0$.

Using both results, necessary and sufficient conditions are given for that problem and the removing step is avoided.
\end{remark}

\begin{remark}[Single-Objective Case]
The same approach can be applied to solve single-objective
problems. In this case, KKT (or FJ) conditions can be applied
directly to the original problem, without scalarizations.
\end{remark}
\section{Obtaining nondominated solutions by multiobjective
optimality conditions} \label{sec:alg4}

In this section, we address the solution of \ref{eq:mob-nonlinear} by directly applying necessary conditions for multiobjective
problems. With these conditions we do not need to scalarize the
problem, as in the above section, avoiding some steps in the process followed in the previous sections.

The following result states the Fritz-John necessary optimality
conditions for multiobjective problems.

\begin{theo}[Multiobjective FJ necessary conditions, Theorem 3.1.1. in \cite{miettinen99}]
\label{mofj_cond} Consider the problem:
\begin{equation}
\label{eq:nl-mofj}
\begin{array}{lrl}
\min & (f_1(x), \ldots, f_k(x)) & \\
s.t.& g_j(x) &\leq 0 \quad  =1, \ldots, m\\
 & h_r(x) &= 0 \quad r=1, \ldots, s\\
 & x &\in \R^n
 \end{array}
 \end{equation}

Let $x^*$ a feasible solution. Suppose that $f_i$, for $i=1, \ldots, k$, $g_j$, for $j=1,
\ldots, m$ and $h_r$, for $r=1, \ldots, s$, are continuously
differentiable at $x^*$. If $x^*$ is a nondominated solution for
Problem \ref{eq:nl-mofj}, then there exist scalars $\nu_i$, for
$i=1, \ldots, k$, $\lambda_j$, for $j=1, \ldots, m$, and $\mu_r$,
for $r=1, \ldots, s$, such that

\begin{equation}
\label{eq:mofj} \tag{MO-FJ}
\begin{array}{rll}
\dsum_{i=1}^k \nu_i \nabla f_i(x^*) + \dsum_{j=1}^m
\lambda_j\,\nabla g_j(x^*)\; +& \dsum_{r=1}^s \mu_r\,\nabla
h_r(x^*)
&=0\\
\lambda_j\,g_j(x^*)\;=&0 & \mbox{for $j=1, \ldots, m$}\\
\lambda_j \;\geq& 0 & \mbox{for $j=1, \ldots, m$}\\
\nu_i \;\geq& 0 & \mbox{for $i=1, \ldots, k$}\\
(\mathbf{\nu}, \mathbf{\lambda}, \mathbf{\mu}) &\neq (\mathbf{0},
\mathbf{0}, \mathbf{0})
\end{array}
\end{equation}
\end{theo}

With this result, one can solve the system given by the necessary
conditions to obtain candidates to be nondominated solutions for
our problem. For solving this system, we use lexicographical
Gröbner bases as in the above sections. We summarize the algorithm
for solving the multiobjective polynomial problem in Algorithm
\ref{alg4}.

\begin{algorithm}[h]
\label{alg4} \SetLine \SetKwInOut{Input}{Input}
\SetKwInOut{Output}{Output}

\Input{$f_1, \ldots, f_k$, $g_1, \ldots g_m, h_1, \ldots, h_r \in \R[x_1, \ldots,
x_n]$}

\textbf{Algorithm: }
\begin{description}
\item[Step 1] Solve system \eqref{eq:mofj}: $X_E^{MOFJ}$.
\item[Step 2] Remove from $X_E^{MOFJ}$ the subset of dominated solutions: $X_E$.
\end{description}

\Output{$X_E$ the set of nondominated solutions for Problem
\ref{eq:mob-nonlinear}}

 \caption{Summary of the procedure for solving MOPBP using the multiobjective FJ optimality conditions.}
\end{algorithm}

\begin{remark}
In the special case where both objective functions and constraints are convex, Theorem \ref{mofj_cond} gives sufficient nondominance conditions for \ref{eq:mob-nonlinear} requiring that $\nu_i>0$ (see Theorem 3.1.8 in \cite{miettinen99}).
\end{remark}

\section{Computational Experiments} \label{sec:comp}

A series of computational experiments have been performed in order
to evaluate the behavior of the proposed solution methods.
Programs have been coded in MAPLE 11 and executed in a PC with an
Intel Core 2 Quad processor at 2x 2.50 Ghz and 4 GB of RAM. The
implementation has been done in that symbolic programming language,
available upon request, in order to make the access easy to both
optimizers and algebraic geometers.

We run the algorithms for three families of binary biobjective and triobjective knapsack problems: linear, quadratic and cubic, and for a biobjective portfolio selection model. For each problem, we obtain the set of nondominated solutions as well as the CPU times for computing the corresponding Gröbner bases associated to the problems, and the total CPU times for obtaining the set of solutions.

We give a short description of the problems where we test the algorithms. In all cases, we use binary variables $x_j$, $j=1, \ldots, n$, where $x_j=1$ means that the item (resp. security) $j$ is selected for the knapsack (resp. portfolio) problem.

\begin{enumerate}
\item {\it Biobjective (linear) knapsack problem} (\texttt{biobj\underline{ }linkn}): Assume that $n$ items are given. Item $j$ has associated costs $q_j^1$, $q_j^2$ for two different targets, and a unit profit $a_j$, $j=1, \ldots, n$. The biobjective knapsack problem calls for selecting the item subsets whose overall profit ensures a knapsack with value at least $b$, so as to minimize (in the nondominance sense) the overall costs. The problem may be formulated:
\begin{eqnarray*}
\min \; (\dsum_{j=1}^n q_{j}^1\,x_j, \dsum_{j=1}^n q_{j}^2\,x_j)\nonumber\\
s.t. \quad\dsum_{i=1}^n\,a_{i}\,x_i \geq b, &  x \in \{0,1\}^n\nonumber
 \end{eqnarray*}
 \item {\it Biobjective cubic knapsack problem} (\texttt{biobj\underline{ }cubkn}): Assume that $n$ items are given where item $j$ has an integer profit $a_j$. In addition we are given two $n\times n \times n$ matrices $P^1 = (p^1_{ijk})$ and $P^2 = (p^2_{ijk})$, where $p^1_{ijk}$ and $p^2_{ijk}$ are the costs for each of the targets if the combination of items $i, j, k$ is selected for $i<j<k$; and two additional $n\times n$ matrices $Q^1 = (q^1_{ij})$ and $Q^2 = (q^2_{ij})$, where $q^1_{jj}$ and $q_{ij}^2$ are the costs for the two different targets if both items $i$ and $j$ are selected for $i<j$  The biobjective cubic knapsack problem calls for selecting the item subsets whose overall profit exceeds the purpose of the knapsack $b$, so as to minimize the overall costs. The problem may be formulated:
\begin{eqnarray*}
\min \; \big(\dsum_{i=1}^{n}\dsum_{j=i}^n q_{ij}^1\,x_i\,x_j + \dsum_{i=1}^{n-2}\dsum_{j=i+1}^{n-1}\dsum_{l=j+1}^n\,p_{ijl}^1\,x_i\,x_j\,x_l, &
 \dsum_{i=1}^{n}\dsum_{j=i}^n q_{ij}^2\,x_i\,x_j + \dsum_{i=1}^{n-2}\dsum_{j=i+1}^{n-1}\dsum_{l=j+1}^n\,p_{ijl}^2\,x_i\,x_j\,x_l\big)\\
s.t. \quad \dsum_{i=1}^n\,a_{i}\,x_i \geq b, \quad x \in \{0,1\}^n
 \end{eqnarray*}
 \item {\it Biobjective quadratic knapsack problem} (\texttt{biobj\underline{ }qkn}): This problem may be seen as a special case of the biobjective cubic knapsack problem when there are no cost correlations between triplets.
\item {\it Triobjective (linear) knapsack problem} (\texttt{triobj\underline{ }linkn}): Assume that $n$ items are given where item $j$ has an integer profit $a_j$. In addition, we are given three vectors $q^1 = (q^1_{j})$, $q^2 = (q^2_{j})$ and  $q^3 = (q^3_{j})$, where $q^1_{j}$,$q_{j}^2$ and $q_{j}^3$ are the costs for three different targets if $j$ is selected. The triobjective knapsack problem calls for selecting the item subsets whose overall profit ensures a profit for the knapsack at least $b$, so as to minimize (in the nondominance sense) the overall costs. The problem is:
\begin{eqnarray*}
\min \; (\dsum_{j=1}^n q_{j}^1\,x_j, \dsum_{j=1}^n q_{j}^2\,x_j, \dsum_{j=1}^n q_{j}^3\,x_j)\\
\qquad s.t.\dsum_{i=1}^n\,a_{i}\,x_i \geq b, \quad x \in \{0,1\}^n
 \end{eqnarray*}
\item {\it Triobjective cubic knapsack problem} (\texttt{triobj\underline{ }cubkn}): We are given $n$ items where item $j$ has an integer profit $a_j$. In addition, we are given three $n\times n \times n$ matrices $P^1 = (p^1_{ijk}), P^2 = (p^2_{ijk})$ and $P^3 = (p^3_{ijk})$, where $p^1_{ijk}, p^2_{ijk}$ and $p^3_{ijk}$ are the costs for each of the targets if the combination of items $i, j$ and $k$ is selected for $i<j<k$; and three additional $n\times n$ matrices $Q^1 = (q^1_{ij}), Q^2 = (q^2_{ij})$ and $Q^3 = (q^2_{ij})$, where $q^1_{jj}, q_{ij}^2$ and $q_{ij}^3$ are the costs for three different targets if both items $i$ and $j$ are selected for $i<j$  The triobjective cubic knapsack problem calls for selecting the item subsets whose overall profit ensures a value of $b$, so as to minimize the overall costs. The problem is:
{\begin{eqnarray*}
\min \; &\big(\dsum_{i=1}^{n}\dsum_{j=i}^n q_{ij}^1\,x_i\,x_j + \dsum_{i=1}^{n-2}\dsum_{j=i+1}^{n-1}\dsum_{l=j+1}^n\,p_{ijl}^1\,x_i\,x_j\,x_l, \dsum_{i=1}^{n}\dsum_{j=i}^n q_{ij}^2\,x_i\,x_j + \dsum_{i=1}^{n-2}\dsum_{j=i+1}^{n-1}\dsum_{l=j+1}^n\,p_{ijl}^2\,x_i\,x_j\,x_l,\\ &\dsum_{i=1}^{n}\dsum_{j=i}^n q_{ij}^3\,x_i\,x_j + \dsum_{i=1}^{n-2}\dsum_{j=i+1}^{n-1}\dsum_{l=j+1}^n\,p_{ijl}^3\,x_i\,x_j\,x_l \big)\\
&s.t.\dsum_{i=1}^n\,a_{i}\,x_i \geq b, \quad x \in \{0,1\}^n
 \end{eqnarray*}}
\item {\it Triobjective quadratic knapsack problem} (\texttt{triobj\underline{ }qkn}): This problem may be seen as a special case of the triobjective cubic knapsack problem when there are no cost correlations between triplets.
\item {\it Biobjective portfolio selection} (\texttt{portfolio}): Consider a market with $n$ securities. An investor with initial wealth $b$ seeks to improve his wealth status by investing it into these $n$ risky securities. Let $X_i$ be the random return per a lot of the $i$-th secutiry ($i=1, \ldots, n$). The mean, $\mu_i = E[X_i]$, and the covariance, $\sigma_{ij}= Cov(X_i, X_j), i, j=1, \ldots, n$, of the returns are assumed to be known. Let $x_i$ be a decision variable that takes value $1$ if the decision-maker invests in the $i$-th security and 0 otherwise. Denote the decision vector by $x=(x_1, \ldots, x_n)$. Then, the random return for a inversion vector $x$ from the securities is $\sum_{i=1}\,x_iX_i$ and the mean and variance of this random variable are $E[\sum_{i=1}\,x_iX_i] = \dsum_{i=1}^n \mu_i x_i$ and $Var(\sum_{i=1}\,x_iX_i) = \dsum_{i=1}^n\dsum_{j=1}^n x_i\,x_j\,\sigma_{ij}$.

Let $a_i$ be the current price of the $i$-th security. Then, if an investor looks for minimizing his investment risk and simultaneously maximizing the expected return with that investment, the problem can be formulated as:
\begin{center}
\begin{tabular}{lcl}
\begin{minipage}[h]{6cm}
$\left.\begin{array}{l}
\min \; ( Var(\dsum_{i=1}^n x_i\,X_i), -E[\dsum_{i=1}^n x_i\,X_i] )\\
s.t.\dsum_{i=1}^n\,a_{i}\,x_i \leq b, \quad x \in \{0,1\}^n
\end{array}\right\}$
\end{minipage}
& $\Rightarrow$ &
\begin{minipage}[h]{6cm}
$\left.\begin{array}{l}
\min \; ( \dsum_{i=1}^n \sigma_{ij}\,x_i\,x_j, -\dsum_{i=1}^n \mu_i\,x_i )\\
 s.t.\dsum_{i=1}^n\,a_{i}\,x_i \leq b, \quad x \in \{0,1\}^n
\end{array}\right\}$
\end{minipage}
\end{tabular}
\end{center}
\end{enumerate}
For each of the above 7 classes of problems, we consider instances randomly generated as follows:  $a_i$ is randomly drawn in $[-10, 10]$
and the coefficients of the objective functions, $q_{ij}^k$, $p^k_{ijl}$, $\sigma_{ij}$ and $\mu_i$, range in $[-10, 10]$.
Once the constraint vector, $(a_1, \ldots, a_n)$, is generated, the right hand side, $b$, is randomly generated in $[1, |\sum_{i=1}^n a_i|]$. For each type of instances and each value of $n$ in $[2, 13]$ we generated 5 instances.

Tables \ref{comp:tabla1} and \ref{comp:tabla2} contain a summary of the average results over the different instances generated for the above problems. Each algorithm is labeled conveniently: \texttt{alg1} corresponds with Algorithm \ref{alg1}, \texttt{kkt} is Algorithm \ref{alg:kkt}, \texttt{kkt\underline{ }sl} is Algorithm \ref{alg:kkt} where the inequality is transformed to an equation using a slack variable, \texttt{fj} is Algorithm \ref{alg:fj}, \texttt{fj\underline{ }sl} is Algorithm \ref{alg:fj} where the inequality is transformed to an equation using a slack variable and \texttt{mofj} stands for Algorithm \ref{alg4}. For each of these algorithms we present the CPU time for computing the corresponding Gröbner basis (\texttt{tgb}), the total CPU time for obtaining the set of nondominated solutions (\texttt{ttot}), the number of nondominated solutions (\texttt{\#nd}) and the number of variables involved in the resolution of the problem (\texttt{\#vars}).

From those tables, the reader may note that Algorithm \ref{alg1} is faster than the others for the smallest instances, although the CPU times for this algorithm increase faster than for the others and it is not able to obtain solutions when the size of the problem is around 12 variables. The algorithms based on Chebyshev scalarization (\texttt{kkt}, \texttt{kkt\underline{ }sl}, \texttt{fj} and \texttt{fj\underline{ }sl}) are better than \texttt{alg1} for the largest instances. The differences between these four methods are meaningful, but the algorithms based on the KKT conditions are, in almost all the instances, faster than those based on the FJ conditions. Note that considering slack variables to avoid the inequality constraint is not better, since the CPU times when the slack variable is considered are larger. Finally, the best algorithm, in CPU time, is \texttt{mofj} since except for the small instances is the fastest and it was able to solve larger instances.

One may think that the last step of our methods, i.e. removing dominated
solutions, should be more time consuming in \texttt{alg1} than in the remaining
methods since \texttt{alg1} does not use optimality conditions. However, from
our experiments this conclusion is not clearly supported. Actually, although this
process, in time consuming, when the dimension of the problem increases this
time is rather small compared with the effort necessary to obtain the Gröbner
bases.

\begin{sidewaystable}[t]
\vspace*{15cm}
 {\scriptsize
\begin{center}
\begin{tabular}{|c|l|lll|lll|lll|lll|lll|lll|l|}\hline
& & \multicolumn{3}{|c|}{\texttt{alg1}} &  \multicolumn{3}{|c|}{\texttt{kkt}} & \multicolumn{3}{|c|}{\texttt{kkt\underline{ }sl}} &\multicolumn{3}{|c|}{\texttt{fj}} & \multicolumn{3}{|c|}{\texttt{fj\underline{ }sl}} & \multicolumn{3}{|c|}{\texttt{mofj}} & \\\hline
prob  &\texttt{n} &  \texttt{gbt} & \texttt{tott} & \texttt{\#vars} & \texttt{gbt} & \texttt{tott} & \texttt{\#vars} & \texttt{gbt} & \texttt{tott} & \texttt{\#vars} & \texttt{gbt} & \texttt{tott} & \texttt{\#vars} & \texttt{gbt} & \texttt{tott} & \texttt{\#vars} & \texttt{gbt} & \texttt{tott} & \texttt{\#vars} &\texttt{\#nd}\\\hline
\multirow{11}{*}{\begin{sideways}\texttt{biobj\underline{ }linkn}\end{sideways}} &2 & 0.02 & 0.07 &  5 & 0.30 & 0.48 & 11 & 0.35 & 0.61 & 12 & 0.24 & 0.39 & 10 & 0.31 & 0.50 & 11 & 0.03 & 0.06 &  7 & 1.8 \\
 &3 & 0.03 & 0.09 &  6 & 0.82 & 1.31 & 13 & 0.95 & 1.61 & 14 & 0.71 & 1.06 & 12 & 0.84 & 1.39 & 13 & 0.10 & 0.18 &  9 & 1.6 \\
 &4 & 0.07 & 0.23 &  7 & 2.45 & 3.89 & 15 & 3.05 & 4.90 & 16 & 2.07 & 3.10 & 14 & 2.54 & 4.12 & 15 & 0.33 & 0.52 & 11 & 3.6 \\
 &5 & 0.32 & 0.69 &  8 & 7.01 & 10.82 & 17 & 8.71 & 13.90 & 18 & 5.92 & 8.71 & 16 & 7.31 & 11.71 & 17 & 0.91 & 1.40 & 13 & 4.6 \\
 &6 & 2.86 & 4.11 &  9 & 25.00 & 37.88 & 19 & 32.29 & 50.13 & 20 & 20.63 & 30.11 & 18 & 26.66 & 41.70 & 19 & 3.01 & 4.61 & 15 & 5.4 \\
 &7 & 31.15 & 35.59 & 10 & 72.85 & 107.15 & 21 & 82.90 & 127.08 & 22 & 55.68 & 79.80 & 20 & 66.77 & 104.82 & 21 & 7.89 & 11.60 & 17 &  5 \\
 &8 & 342.10 & 373.72 & 11 & 176.94 & 261.58 & 23 & 209.68 & 322.58 & 24 & 133.78 & 193.54 & 22 & 165.97 & 262.46 & 23 & 18.35 & 27.60 & 19 & 4.2 \\
 &9 & 5273.56 & 6382.43 & 12 & 462.14 & 675.24 & 25 & 529.50 & 813.54 & 26 & 333.81 & 492.89 & 24 & 418.75 & 670.82 & 25 & 48.50 & 79.08 & 21 &  8 \\
 &10 & & & & & & & & & & & & & & & & 269.19 & 404.04 & 23 & 7.8 \\
 &11 &   &   &   &   &   &   &   &   &   &   &   &   &   &   &   &   480.26 &   835.46 &  25 & 6.4 \\
 &12 &   &   &   &   &   &   &   &   &   &   &   &   &   &   &   &    1340.31 &    2004.33 &  27 & 5.4 \\
 &13 &   &   &   &   &   &   &   &   &   &   &   &   &   &   &   & 4091.92 & 19546.05 &    29 &    11 \\\hline
\multirow{11}{*}{\begin{sideways}\texttt{biobj\underline{ }qkn}\end{sideways}}  &2 & 0.03 & 0.07 & 5 & 0.28 & 0.47 &11 & 0.37 & 0.64 &12 & 0.23 & 0.36 &10 & 0.31 & 0.51 &11 & 0.04 & 0.07 & 7 &  1.4 \\
 &3 & 0.04 & 0.10 & 6 & 0.79 & 1.29 &13 & 1.12 & 1.86 &14 & 0.68 & 1.04 &12 & 0.95 & 1.54 &13 & 0.09 & 0.16 & 9 & 2 \\
 &4 & 0.19 & 0.37 & 7 & 3.27 & 4.97 &15 & 4.63 & 7.15 &16 & 2.69 & 4.01 &14 & 3.89 & 5.95 &15 & 0.49 & 0.70 &11 &  2.4 \\
 &5 & 1.76 & 2.37 & 8 & 9.22 &13.88 &17 &11.94 &18.06 &18 & 7.54 &10.94 &16 & 9.90 &15.11 &17 & 1.08 & 1.74 &13 &  3.8 \\
 &6 &21.43 &23.22 & 9 &25.21 &37.33 &19 &33.74 &50.13 &20 &20.57 &29.46 &18 &27.60 &41.86 &19 & 2.92 & 4.58 &15 &  3.6 \\
 &7 &  425.13 &  430.43 &10 &59.57 &91.29 &21 &85.56 &  129.88 &22 &49.38 &72.81 &20 &69.20 &  107.76 &21 & 5.93 & 9.71 &17 &  4.4 \\
 &8 &  &  &  &172.26 &255.22 & 23.00 &228.66 &337.02 & 24.00 &127.43 &188.32 & 22.00 &174.23 &272.72 & 23.00 & 14.75 & 25.66 & 19.00 &  5.4 \\
 &9 &   &  &   & 463.39 & 692.25 & 25 & 641.29 & 939.99 & 26 & 350.34 & 515.07 &  24 & 489.31 & 755.88 &    25 & 39.90 & 65.78 &    21 &   6.6 \\
 &10 &  &  &  &  &  &  &  &  &  &  &  &  &  &  &  &138.11 &255.42 & 23.00 &  7 \\
 &11 &    &    &    &    &    &    &    &    &    &    &    &    &    &    &    & 331.34 &     643.90 &  25 &  9 \\
 &12 &    &    &    &    &    &    &    &    &    &    &    &    &    &    &    &  891.46 &    1833.75 &  27 &  8 \\\hline
 \multirow{9}{*}{\begin{sideways}\texttt{biobj\underline{ }cubkn}\end{sideways}}  &   3 & 0.17 & 0.21 &    7 & 0.72 & 1.14 &   13 & 0.97 & 1.56 &   14 & 0.62 & 0.92 &   12 & 0.82 & 1.31 &   13 & 0.09 & 0.15 &    9 &  1.8 \\
 &  4 & 1.90 & 2.00 &    8 & 4.68 & 6.73 &   15 & 6.26 & 8.98 &   16 & 3.86 & 5.37 &   14 & 5.12 & 7.38 &   15 & 0.56 & 0.83 &   11 &  1.8 \\
 & 5 &  7.26 &  7.48 &  9 &  9.29 & 13.43 & 17 & 12.49 & 17.79 & 18 &  7.86 & 10.81 & 16 & 10.31 & 14.91 & 17 &  1.11 &  1.65 & 13 &  2.8 \\
 &   6 &     205.52 &     206.25 &   10 &&&&&&&&&&&&& 6.11 & 8.29 &   15 &  2.2 \\
 &   7 &    2067.24 &    2069.38 &   11 &&&&&&&&&&&&&11.25 &15.24 &   17 &  6.6 \\
 &   8 &&&&&&&&&&&&&&&&24.98 &32.25 &   19 &  5.2 \\
 &   9 &&&&&&&&&&&&&&&&81.80 &     106.90 &   21 &  6.2 \\
 &   10 &   &   &   &   &   &   &   &   &   &   &   &   &   &   &   &  258.68 &     389.95 &23 & 5.8 \\
 &   11 &   &   &   &   &   &   &   &   &   &   &   &   &   &   &   &  690.43 &     916.80 &25 & 7.8 \\\hline
 \end{tabular}
 \caption{\label{comp:tabla1}Computational results for biobjective knapsack problems.}
  \end{center}}
\begin{center}
\begin{tabular}{|l|c|c|c|}\hline
Algorithm & \texttt{\#var} & \texttt{\#gen} & \texttt{maxdeg}\\\hline
\texttt{alg1} & $2n+k+m+s$ & $n+k+m+s$             &$\max\{2, deg(f), deg(g), deg(h)\}$            \\\hline
\texttt{kkt}  & $2n+2k+m+s+1$      & $2n+k+m+s+1$   &       $\max\{deg(f)+2, deg(g)+1, deg(h)\}$            \\
\texttt{nr}   & $2n+2k+m+s+1$       & $2n+m+s$       &  $\max\{ deg(f)+1, deg(g), deg(h)\}$            \\\hline
\texttt{fj}            & $2n+2k+m+s+2$       &  $2n+k+m+s+1$  & $\max\{deg(f)+2, deg(g)+1, deg(h)\}$            \\\hline
\texttt{mojf}        & $2n+k+m+s$       & $2n+m+s$       & $\max\{deg(f),deg(g)+1, deg(h)\}$            \\\hline
\end{tabular}
\caption{\label{tabla3}Information about all the algorithms.}
  \end{center}
 \end{sidewaystable}
\begin{sidewaystable}[htb!]
\vspace*{16cm}
{\scriptsize
\begin{center}
\begin{tabular}{|c|l|lll|lll|lll|lll|lll|lll|l|}\hline
& & \multicolumn{3}{|c|}{\texttt{alg1}} &  \multicolumn{3}{|c|}{\texttt{kkt}} & \multicolumn{3}{|c|}{\texttt{kkt\underline{ }sl}} &\multicolumn{3}{|c|}{\texttt{fj}} & \multicolumn{3}{|c|}{\texttt{fj\underline{ }sl}} & \multicolumn{3}{|c|}{\texttt{mofj}} & \\\hline
prob  &\texttt{n} &  \texttt{gbt} & \texttt{tott} & \texttt{\#vars} & \texttt{gbt} & \texttt{tott} & \texttt{\#vars} & \texttt{gbt} & \texttt{tott} & \texttt{\#vars} & \texttt{gbt} & \texttt{tott} & \texttt{\#vars} & \texttt{gbt} & \texttt{tott} & \texttt{\#vars} & \texttt{gbt} & \texttt{tott} & \texttt{\#vars} &\texttt{\#nd}\\\hline
\multirow{11}{*}{\begin{sideways}\texttt{triobj\underline{ }linkn}\end{sideways}}  &  2 &0.03 &0.08 & 6 &0.82 &1.27 &13 &0.99 &1.61 &14 &0.75 &1.16 &12 &0.88 &1.47 & 13 &0.04 &0.06 & 8 & 1.6 \\
 &  3 &0.02 &0.10 & 7 &2.98 &4.33 &15 &3.63 &5.51 &16 &4.19 &5.62 &14 &4.95 &7.00 & 15 &0.14 &0.20 &10 & 2 \\
 &  4 &0.06 &0.65 & 8 & 14.71 & 20.07 &17 & 18.03 & 25.22 &18 & 19.93 & 25.61 &16 & 23.62 & 31.94 & 17 & 0.73 & 0.90 &12 & 2.8 \\
 & 5 &0.71 &1.00 & 9 & 40.75 & 56.40 &19 & 48.20 & 68.21 &20 & 55.31 & 72.04 &18 & 65.38 & 87.07 & 19 & 2.70 & 3.62 &14 & 4.2 \\
 & 6 &2.44 &3.32 &10 & 68.69 & 102.22 &21 & 84.95 & 129.12 &22 & 107.04 & 147.24 &20 & 130.13 & 185.43 & 21 & 2.26 & 3.36 &16 & 6.2 \\
 &  7 & 22.35 & 25.32 &11 & 188.86 & 276.15 &23 & 219.78 & 335.08 &24 & 272.32 & 374.25 &22 & 326.76 & 470.51 & 23 &5.89 &8.59 &18 &11.4 \\
 & 8 & 360.38 & 367.76 &12 & 501.35 & 731.46 &25 & 576.07 & 856.28 &26 & 729.93 & 993.07 &24 & 830.08 &    1194.75 & 25 & 17.16 & 23.83 &20 & 4.4 \\
 & 9 &   &   &   &   &   &   &   &   &   &   &   &   &   &   &   & 69.68 & 115.30 &22 &24.4 \\
 & 10 &   &   &   &   &   &   &   &   &   &   &   &   &   &   &   & 193.15 & 301.04 &24 &31.4 \\
 & 11 &    &    &    &    &    &    &    &    &    &    &    &    &    &    &    &     529.51 &    1075.93 & 26 & 26 \\
 & 12 &   &   &   &   &   &   &   &   &   &   &   &   &   &   &   &    1422.70 &  3145.30 &28 &  85.2 \\
\hline
\multirow{11}{*}{\begin{sideways}\texttt{triobj\underline{ }qkn}\end{sideways}}  & 2 & 0.04 & 0.09 &  6 & 0.96 & 1.52 & 13 & 1.28 & 1.96 & 14 & 0.88 & 1.37 & 12 & 1.19 & 1.81 & 13 & 0.05 & 0.08 &  8 &  2.2 \\
 & 3 & 0.06 & 0.18 &  7 & 3.20 & 4.72 & 15 & 4.11 & 6.13 & 16 & 4.67 & 6.25 & 14 & 5.21 & 7.42 & 15 & 0.14 & 0.20 & 10 &  2.8 \\
 & 4 & 0.09 & 0.32 &  8 & 6.50 & 9.89 & 17 & 9.03 &  13.82 & 18 & 9.59 &  13.14 & 16 &  11.44 &  16.59 & 17 & 0.28 & 0.45 & 12 &  5.2 \\
 & 5 & 2.87 & 3.68 &  9 &  21.05 &  30.87 & 19 &  30.61 &  45.61 & 20 &  30.93 &  41.79 & 18 &  39.81 &  56.92 & 19 & 0.79 & 1.26 & 14 &  8.6 \\
 & 6 &  31.76 & 33.69 & 10 &  49.11 & 72.24 & 21 &  72.03 & 106.16 & 22 &  69.31 &  91.93 & 20 &  86.83 & 122.30 & 21 & 2.19 & 3.85 & 16 &  9.4 \\
 & 7 & 1099.24 & 1109.07 & 11 & 152.63 & 224.47 & 23 & 223.51 & 326.03 & 24 & 232.64 & 319.26 & 22 & 296.38 & 417.00 & 23 & 5.18 & 8.14 & 18 & 12.6 \\
 & 8 & & & & 481.73 & 709.09 & 25 & 721.64 & 1078.26 & 26 & 700.63 & 980.21 & 24 & 904.31 & 1289.78 & 25 & 14.81 & 22.22 & 20 & 17.6 \\
 & 9 &  &  &  &  &  &  &  &  &  &  &  &  &  &  &  &  44.36 &  67.14 & 22 & 13.6 \\
 & 10 &  &  &  &  &  &  &  &  &  &  &  &  &  &  &  & 119.50 & 213.66 & 24 & 17.8 \\
 & 11 &   &   &   &   &   &   &   &   &   &   &   &   &   &   &   &     343.53 & 793.88 &26 &  29.8 \\
 & 12 &   &   &   &   &   &   &   &   &   &   &   &   &   &   &   &  1018.35 &  2445.64 & 28 & 40.2 \\\hline
\multirow{9}{*}{\begin{sideways}\texttt{triobj\underline{ }cubkn}\end{sideways}}& 3 & 0.05 & 0.15 & 7 & 2.74 & 4.01 & 15 & 4.07 & 5.88 & 16 & 4.25 & 5.55 & 14 & 5.33 & 7.29 & 15 & 0.11 & 0.19 & 10 &  3.4 \\
& 4 & 0.18 & 0.43 & 8 & 9.97 &  13.81 & 17 &  13.14 &  18.50 & 18 &  14.32 &  18.31 & 16 &  17.91 &  23.90 & 17 & 0.42 & 0.56 & 12 &  4.4 \\
& 5 & 1.37 & 1.82 & 9 &  &  & &  &  &  &  &  &  &  &  & & 1.03 & 1.43 & 14 &  8.4 \\
& 6 &  28.61 &  29.80 & 10 &  &  &  &  &  & &  &  &  &  &  &  & 2.93 & 4.19 & 16 & 6 \\
& 7 &  &  &  &  &  &  &  &  &  &  &  &  &  &  &  & 9.57 &  12.67 & 18 &  7.6 \\
& 8 &  &  &  &  &  &  &  &  &  &  &  &  &  &  &  &  30.95 &  37.29 & 20 & 13.4 \\
& 9 &  &  &  &  &  &  &  &  &  &  &  &  &  &  &  &  92.52 & 115.68 & 22 & 25.4 \\
&10 &   &   &   &   &   &   &   &   &   &   &   &   &   &   &   & 371.98 & 447.80 &24 & 29.4 \\
&11 &&&&&&&&&&&&&&&& 1006.29 & 1593.48 &    26 &  37.4 \\
\hline
\multirow{10}{*}{\begin{sideways}\texttt{portfolio}\end{sideways}}& 2 &  0.03 &  0.14 & 5 &  0.43 &  0.73 &11 &  0.68 &  1.10 &12 &  0.40 &  0.65 &10 &  0.60 &  0.94 &11 &  0.06 &  0.08 & 7 &  2.0 \\
& 3 &  0.03 &  0.08 & 6 &  1.40 &  2.09 &13 &  1.65 &  2.59 &14 &  1.08 &  1.64 &12 &  1.42 &  2.28 &13 &  0.18 &  0.28 & 9 &  2.0\\
& 4 &  0.08 &  0.38 & 7 &  3.50 &  5.34 &15 &  4.56 &  7.08 &16 &  2.85 &  4.17 &14 &  3.74 &  5.95 &15 &  0.46 &  0.71 &11 &  4.0 \\
& 5 &  0.82 &  1.33 & 8 &  9.41 & 14.46 &17 & 11.87 & 18.16 &18 &  7.85 & 11.61 &16 &  9.76 & 15.15 &17 &  1.16 &  1.84 &13 &  5.4 \\
& 6 & 17.39 & 18.73 & 9 & 25.28 & 38.45 &19 & 32.75 & 50.01 &20 & 20.92 & 30.43 &18 & 26.76 & 41.96 &19 &  2.90 &  4.71 &15 &  6.8 \\
& 7 & 179.84 & 184.29 &10 & 74.90 & 111.93 &21 & 95.74 & 143.34 &22 & 59.19 & 85.91 &20 & 75.66 & 117.77 &21 &  7.40 & 12.09 &17 &  6.2 \\
& 8 & 4739.76 & 4749.91 &11 & 132.97 & 201.62 &23 & 168.92 & 255.92 &24 & 108.75 & 161.05 &22 & 137.98 & 217.42 &23 & 15.16 & 25.23 &19 &  8.0 \\
&9 &  &  &  & 393.58 & 606.13 & 25 & 580.04 & 886.54 & 26 & 311.33 & 490.41 & 24 & 455.34 & 742.80 & 25 &  35.83 &  76.77 & 21 & 11.4 \\
&10 &  &  &  &  1212.90 &  1837.36 & 27 & 1601.12 & 2440.86 & 28 & 869.72 & 1379.88 & 26 &  1207.49 &  1982.07 & 27 & 101.17 & 221.92 & 23 & 12.2 \\
&  11 &  &  &  &  &  &  &  &  &  &  &  &  &  &  &  & 305.80 & 724.28 & 25 & 15.6 \\
\hline
 \end{tabular}
 \end{center}
 \caption{\label{comp:tabla2}Computational results for triobjective knapsack and biobjective portfolio problems.}}
 \end{sidewaystable}

Table \ref{tabla3} shows some information about each of the presented algorithms. For a multiobjective problem with $n$ variables, $m$ polynomial inequality constraints given by $\mathbf{g} = (g_1, \ldots, g_m)$, $s$ polynomial equality constraints given by $\mathbf{h} = (h_1, \ldots, h_s)$ and $k$ objectives functions given by $\mathbf{f}=(f_1, \ldots, f_k)$, Table 3 shows the number of variables (\texttt{\#var}), the number of generators (\texttt{\#gen}) and the maximal degrees (\texttt{maxdeg}) of the initial polynomial ideals related to the each of the algorithms. These numbers inform us about the theoretical complexity of the algorithms. The computation of Gröbner bases depends of the number of variables (in general, double exponential) and of the size of the initial system of generators (degrees and number of polynomials). Actually, it is known that computing a Gröbner basis using Buchberger Algorithm is doubly exponential in the number of variables. Some complexity bounds for this algorithm involving  \texttt{\#var}, \texttt{\#gen} and \texttt{maxdeg} can be found in \cite{dube-mishra-yap95}.

From the above table the  reader may note that both \texttt{alg1} and \texttt{mojf}  have the
same number of variables in any case, but the number of initial generators for \texttt{alg1} is,
in general, smaller than the same number for \texttt{mojf}, since the number of objectives is usually smaller than the number of variables. Furthermore, maximal degrees are smaller in \texttt{alg1} than in \texttt{mojf} . However, in practice, \texttt{mojf}  is faster than \texttt{alg1} since using optimality conditions helps in identifying nondominated solutions.

\end{document}